\documentclass[11pt,twoside,english]{myamsart}
\advance\oddsidemargin by -1.0cm
\advance\evensidemargin by -1.0cm
\advance\topmargin by -1.0cm 
\usepackage{amssymb}
\usepackage{babel}
\usepackage{amstext}
\usepackage{amscd}   
\usepackage{epsfig}  
\usepackage{rotating}
\let\mathg\mathfrak
\theoremstyle{plain}
\newtheorem{cor}{Corollary}[section]

\newtheorem{thm}{Theorem}[section]
\newtheorem{prop}{Proposition}[section]
\theoremstyle{definition}

\newtheorem{NB}{Remark}[section]

\renewcommand{\labelenumi}{$(\theenumi)$}
\newcommand{\bdm}{\begin{displaymath}}
\newcommand{\edm}{\end{displaymath}}
\newcommand{\be}{\begin{equation}}
\newcommand{\ee}{\end{equation}}
\newcommand{\ba}[1]{\begin{array}{#1}}
\newcommand{\ea}{\end{array}}

\newcommand{\btab}{\begin{tabular}}
\newcommand{\etab}{\end{tabular}}

\newcommand{\x}{\times}
\newcommand{\op}{\oplus}

\newcommand{\ra}{\rightarrow}

\newcommand{\tr}{\ensuremath{\mathrm{tr}}}

\newcommand{\R}{\ensuremath{\mathbb{R}}}

\newcommand{\M}{\ensuremath{\mathcal{M}}}

\newcommand{\T}{\ensuremath{\mathrm{T}}}
\newcommand{\F}{\ensuremath{\mathrm{F}}}
\newcommand{\G}{\ensuremath{\mathrm{G}}}

\newcommand{\eps}{\ensuremath{\varepsilon}} 
    
\newcommand{\vrho}{\ensuremath{\varrho}}

\newcommand{\Ad}{\ensuremath{\mathrm{Ad}\,}}

\newcommand{\diag}{\ensuremath{\mathrm{diag}}}

\newcommand{\su}{\ensuremath{\mathg{su}}}
\newcommand{\SU}{\ensuremath{\mathrm{SU}}}

\newcommand{\so}{\ensuremath{\mathg{so}}}
\newcommand{\SO}{\ensuremath{\mathrm{SO}}}
\newcommand{\Spin}{\ensuremath{\mathrm{Spin}}}
\newcommand{\spin}{\ensuremath{\mathg{spin}}}

\newcommand{\g}{\ensuremath{\mathfrak{g}}}

\newcommand{\m}{\ensuremath{\mathfrak{m}}}

\begin{document}
\def\haken{\mathbin{\hbox to 6pt{%
                 \vrule height0.4pt width5pt depth0pt
                 \kern-.4pt
                 \vrule height6pt width0.4pt depth0pt\hss}}}
    \let \hook\intprod
\setcounter{equation}{0}
%
%
\thispagestyle{empty}
%
\date{\today}
\title[Killing spinors in Supergravity with $4$-fluxes]
{Killing spinors in Supergravity with $4$-fluxes}
%
%
%
\author{Ilka Agricola}
\author{Thomas Friedrich}
\address{\hspace{-5mm} 
{\normalfont\ttfamily agricola@mathematik.hu-berlin.de}\newline
{\normalfont\ttfamily friedric@mathematik.hu-berlin.de}\newline
Institut f\"ur Mathematik \newline
Humboldt-Universit\"at zu Berlin\newline
Sitz: WBC Adlershof\newline
D-10099 Berlin, Germany}
%
\thanks{Supported by the SFB 288 "Differential Geometry
and Quantum Physics" of the DFG and the research project
``Special Geometries in Mathematical Physics'' of the Volks\-wagen Foundation}
\subjclass[2000]{Primary 53 C 25; Secondary 81 T 30}
\begin{abstract}
We study the spinorial Killing equation of supergravity involving
a torsion $3$-form $\T$ as well as a flux $4$-form $\F$. In dimension
seven, we construct explicit families of compact solutions out of
$3$-Sasakian geometries, nearly parallel $\G_2$-geometries and
on the homogeneous Aloff-Wallach space. The constraint $\F \cdot \Psi = 0$
defines a non empty subfamily of solutions. We investigate the
constraint $\T \cdot \Psi = 0$, too, and show that it singles out a
very special choice of numerical parameters in the Killing equation,
which can also be justified geometrically. 
\end{abstract}
\maketitle
\tableofcontents
\pagestyle{headings}
%
%
%
\section{Introduction}\noindent
Supergravity models can be described geometrically by a
tuple $(M^n, \, g, \, \T, \, \F, \, \Psi)$ consisting of
a Riemannian manifold, a $3$-form $\T$, a $4$-form flux 
$\F$ and a spinor field $\Psi$. The link between
these geometric objects is the so called Killing equation 
(see \cite{Duff}, \cite{FigueroaPapadopoulus})
\bdm
\nabla^g_X \Psi \, + \, \frac{1}{4} \cdot (X \haken \T) \cdot
\Psi \, + \, \frac{1}{144} \cdot ( X \haken \F \, - \, 
8 \cdot X \wedge \F) \cdot \Psi \ = \ 0 \, .
\edm
This equation should be satisfied for any tangent vector $X$.
Consequently, we deal with a highly overdetermined system
of first order partial differential equations.
The $3$-form $\T$ has an interpretation as the torsion of a 
linear, metric connection
$\nabla$ with totally skew-symmetric torsion.
$11$-dimensional space-time solutions are interesting
and the models with a maximal number of supersymmetries
have been classified (see \cite{FigueroaPapadopoulus}). 
The Kaluza-Klein reduction of $\M$-theory (see \cite{Ali}, 
\cite{BehrndtJeschek}, \cite{BilalDS} and \cite{WNW}) yields that
dimensions $4 \leq n \leq 8$ are of interest, too. However, 
then additional algebraic constraints occur, for example, an
algebraic coupling between the torsion $3$-form or the flux $4$-form
and the spinor field $\Psi$,
\bdm
\T \cdot \Psi \ = \ 0 \, \quad \mathrm{or} \quad 
\F \cdot \Psi \ = \ 0 \ . 
\edm
The aim of this note is to present a geometric method for
solving the equation under consideration. The main idea of our approach
is easy to explain. We start with a Riemannian manifold admitting
a spinor field $\Psi$ of some special type. For this, there are many 
possibilities. The spinor field may 
be a Riemannian Killing spinor (see \cite{Fri1}) on some irreducible
Einstein space,
\bdm
\nabla^g_X \Psi \ = \ \lambda \cdot X \cdot \Psi \, .
\edm 
The spinor field may be a
K\"ahlerian Killing spinor (see \cite{Kirchberg}) defined on some special 
K\"ahler manifold. In odd dimensions, 
we can start with an
$\eta$-Einstein-Sasakian manifold and its
contact Killing spinor (see \cite{FKim}). On
a reductive space, the spinor field may be 
an invariant spinor of the isotropy representation. In any case, triples 
$(M^n, \, g, \, \Psi)$ of the type we need 
have been studied very intensively in mathematics
since more then $20$ years. In particular, the 
dimensions $n \leq 8$ and the corresponding special geometries play
a crucial role. 
The books \cite{BFGK}, \cite{Fri2} contain the results
as well as the relevant references in detail.
Let us moreover assume
that there exists a ``canonical'' family $(\T_{\omega}, \, \F_{\omega})$ of 
forms on $M^n$ depending on some parameter
$\omega \in \Omega$. We consider the system
\bdm 
\Big(\lambda \cdot X \, + \, \frac{1}{4} \cdot (X \haken \T_{\omega}) \, 
+ \, \frac{1}{144} (X \haken \F_{\omega} \, - \, 8 \cdot
X \wedge \F_{\omega}) \Big) \cdot \Psi \ = \ 0 \,,
\edm
which is a highly overdetermined system of $n \cdot 2^{[n/2]}$ 
algebraic equations in the parameters $\omega \in \Omega$.
Furthermore, we can add the equations $\T_{\omega} \cdot \Psi = 0$ or
$\F_{\omega} \cdot \Psi = 0$. 
It is a matter of fact -- quite surprising to the mathematician -- that 
the solutions often constitute a 
non empty subset of the parameter space. We solve the corresponding
equations with the help of standard math computer programs. In this way, we 
obtain families of $3$- and $4$-forms solving the equation
on the manifold we started with. The first interesting dimension 
in presence of a $4$-flux is seven. 
Therefore, we apply the described method to nearly parallel 
$\G_2$-manifolds, to $3$-Sasakian manifolds
and to the $7$-dimensional Aloff-Wallach space. The $\G_2$-case has
been already investigated in $\M$-theory compactifications 
to dimension four (see \cite{BilalDS}). On a nearly parallel $\G_2$-manifold, 
there is a natural flux form.
The second and third case are more flexible and interesting. 
We will construct families of torsion and flux forms with Killing spinors 
out of a $3$-Sasakian structure. 
The underlying metric of a $3$-Sasakian manifold is still 
Einstein (see \cite{BG}). On the Aloff-Wallach space $N(1,1)$ 
we obtain families of non Einstein metrics
equipped with torsion forms, flux forms and Killing spinors.
The method we use can probably be applied to other dimensions
and special geometries, too.\\

\noindent
Let us consider a slightly more general equation depending on
two real parameters $(p,q) \in \R^2$,
\bdm
\nabla^g_X \Psi \, + \, \frac{1}{4} \cdot (X \haken \T) \cdot
\Psi \, + \, p \cdot ( X \haken \F ) \cdot \Psi \,  + \, 
q \cdot (X \wedge \F) \cdot \Psi \ = \ 0 \, .
\edm
We will construct $7$-dimensional solutions of this field equation 
for {\it any} pair
$(p,q) \in \R^2$ of parameters. 
Obviously, if the flux form is non trivial, then the ratio between
the parameters $p$ and $q$ is important. It turns out that
in any dimension $n$ there is 
one distinguished pair of parameters, namely $4 \, p \, - \, (n-4) \, 
q = 0$. This coupling of the parameters plays a special role for
our solutions. It can be motivated as well by the
observation that in this case the action of the Riemannian Dirac 
operator on the spinor depends only on the torsion form $\T$, but
not on the flux form $\F$.
Remark that this is {\it not} the ratio appearing in the original
equation of supergravity, since we consider positive definite
metrics. The equation in arbitrary dimension 
reads then as
\bdm
\nabla^g_X \Psi \, + \, \frac{1}{4} \cdot (X \haken \T) \cdot
\Psi \, + \, \frac{n-4}{4} \cdot ( X \haken \F ) \cdot \Psi \,  + \, 
(X \wedge \F) \cdot \Psi \ = \ 0 \, .
\edm
In Section $5$, we discuss in more detail the 
$7$-dimensional solutions of this special equation which we constructed.
%
\section{Killing Spinors with $4$-Fluxes on $3$-Sasakian Manifolds}
\label{fam-conn}\noindent
%
\noindent
The structure group of a $3$-Sasakian geometry is 
the subgroup $\SU(2) \subset \mathrm{G}_2 \subset \SO(7)$,
the isotropy group of four spinors in dimension seven. We 
describe the subgroup $\SU(2)$ in such a way that the vectors
$e_1, e_2, e_7 \in \R^7$ are fixed. More precisely, let the Lie algebra
$\su(2)$ be generated by the following $2$-forms in $\R^7$:
\bdm
e_{34} \, + \, e_{56} , \quad 
e_{35} \, - \, e_{46} , \quad e_{36} \, + \, e_{56} .
\edm
The real spin representation $\Delta_7$
splits under the action of $\SU(2)$ into a $4$-dimensional trivial
representation $\Delta_7^0$ and the unique non trivial $4$-dimensional
representation $\Delta_7^{1}$. 
We use the standard realization of the $8$-dimensional
$\Spin(7)$-representation as given in \cite[p.97]{BFGK} or 
\cite[p.13]{Fri2} (see \cite{AgFr}, too). Denote by $\Psi_1, \ldots ,
\Psi_8$ its basis. The space $\Delta_7^0$
is spanned by the spinors $\Psi_3, \Psi_4, \Psi_5, \Psi_6$. We consider
the following $\SU(2)$-invariant $2$-forms on $\R^7$:
\bdm
de_1 \ := \ e_{35} \, + \, e_{46}, \quad 
de_2 \ := \ e_{45} \, - \, e_{36} , \quad de_7 \ := \ e_{34} \, - \, e_{56} .
\edm
Using this notation, we introduce a family of invariant $3$-forms in 
$\R^7$ depending on $10$ parameters,
\bdm
\T \ = \ \sum_{i,j = 1,2,7} t_{ij} \cdot e_i \wedge d e_j \, + \, 
t \cdot e_1 \wedge e_2 \wedge e_7.
\edm
The space of $\SU(2)$-invariant $4$-forms on $\R^7$ has also dimension ten,
\bdm
\F \ = \ \sum_{i,j,k=1,2,7} f_{ijk} \cdot e_i \wedge e_j 
\wedge d e_k \, + \, f \cdot e_3 \wedge e_4 \wedge e_5 \wedge e_6 .
\edm
All together, on a $3$-Sasakian manifold there exists a canonical family
$\Omega$ of forms depending on $20$ parameters.
The key point of our considerations is the following 
algebraic observation.
\begin{prop}\label{3Sas}  
Fix  parameters $(p , q) \in \R^2$.
For any spinor $0 \neq \Psi \in \Delta_7^0$, the set of all pairs
$(\T, \, \F)$ consisting of $\SU(2)$-invariant forms and satisfying,
for any vector $X \in \R^7$, the equation
\bdm
\Big(\frac{1}{2}\cdot X \, + \, \frac{1}{4} \cdot (X \haken \T) \, 
+ \, p \cdot (X \haken \F) \, + \, q \cdot 
( X \wedge \F) \Big) \cdot \Psi \ = \ 0 
\edm
is a $7$-dimensional affine space. The condition $\F \cdot \Psi = 0$ defines 
a $6$-dimensional subspace. If $4 \, p \, - \, 3 \, q \neq 0$, 
the condition $\T \cdot \Psi = 0$ defines a $6$-dimensional affine
subspace, too. If $4 \, p \, - \, 3 \, q = 0$, we have $\T \cdot \Psi = 
(14/3) \cdot \Psi$ for all torsion forms in the family. 
Both constraints $\T \cdot \Psi = 0$ and $\F \cdot \Psi = 0$ cannot
be fulfilled simultaneously.
\end{prop} 
\begin{proof} 
Given a spinor $\Psi = a \, \Psi_3 + b \, \Psi_4 + 
c \, \Psi_5 + d \, \Psi_6$, we solve the overdetermined system
with respect to the coefficients
of the $3$-form $\T$ and the coefficients of the $4$-form $\F$. 
It turns out that solutions exist and can be given explicitely. 
For simplicity, we provide the formulas
in case that $\Psi = \Psi_3$ is one of the basic spinors. The set
of all solutions is an affine space parameterized by $f_{127}, \, 
f_{172}, \, f_{177}, \, f_{271}, \, f_{272}, \, f_{277}, \, f$.
The other coordinates are given by the formulas
\begin{eqnarray*}
t_{12} & = &  - \, 4\, (p\, + \, q) \, f_{272} \ = \ - \, t_{21} , \quad 
t_{17} \ = \ - \, 4 \, (p\, + \, q) \, f_{277} \ = \ - \, t_{71} , \\
t_{27} & = & 4 \, (p \, + \, q) \, f_{177} \ = \ t_{72} , \\
t_{11} & = & - \, \frac{2}{3} \, + \, \frac{4 p}{3} \, f \, - \, 
\frac{8 p}{3} \, f_{127} \, + \, \frac{8 p}{3} \, f_{172} \, - \, 
\frac{4}{3}(p \, + \,  3 q) \, f_{271} , \\
t_{22} & = & \frac{2}{3} \, - \, \frac{4 p}{3} \, f \, + \, 
\frac{8 p}{3} \, f_{127} \, - \, \frac{8 p}{3} \, f_{271} \, + \, 
\frac{4}{3}(p \, + \,  3 q) \, f_{172} , \\
t_{77} & = & \frac{2}{3} \, - \, \frac{4 p}{3} \, f \, - \, 
\frac{8 p}{3} \, f_{172} \, - \, \frac{8 p}{3} \, f_{271} \, - \, 
\frac{4}{3}(p \, + \, 3 q) \, f_{127} , \\
t & = & \frac{2}{3} \, + \, \frac{8 p}{3} \, f_{127} \, - \, 
\frac{8 p}{3} \, f_{172} \, - \, \frac{8 p}{3} \, f_{271} \, + \, 
\frac{4}{3}(2 p \, + \, 3 q) \, f_ , \\
f_{122} & = & - \, f_{177} , \quad 
f_{121} \ = \ - \, f_{277} , \quad 
f_{171} \ = \ f_{272} \ .
\end{eqnarray*}
The equation $\F \cdot \Psi_3 = 0$ describes a $6$-dimensional
linear subspace,
\bdm
f \, - 2 \, f_{127} \ + \, 2 \, f_{172} \, + \, 2 \, f_{271} \ = \ 0 .
\edm 
On the other side, the equation $\T \cdot \Psi_3 = 0$ defines a $6$-dimensional
affine subspace
\bdm
- \, 7 \, + \, (8 \, p \, - 6 \, q) \big(f \, - 2 \, f_{127} \ + 
\, 2 \, f_{172} \, + \, 2 \, f_{271}\big) \ = \ 0 .
\edm
The intersection of these spaces is empty.
\end{proof}
\begin{NB} 
The cases $\T = 0$ and $\F = 0$ have been discussed already
in \cite{AgFr}. Under this constraint the algebraic system has a 
unique solution.
\end{NB}
\begin{NB} 
From the geometric point of view, there is an interesting 
case, namely $ 4 \, p \, - 3 \, q = 0$. This is {\it not} the ratio
of the parameters $p,q$ appearing in supergravity. 
In this case, the constraint
$\T \cdot \Psi = 0$ is never satisfied for a non trivial spinor. 
In Section $5$, we will discuss this family of solutions in more detail.
\end{NB} 
 
\noindent
Consider a simply connected $3$-Sasakian manifold $M^7$ of dimension seven and 
denote its three contact structures by $\eta_1, \eta_2$, and $\eta_7$. 
It is known that $M^7$ is then an Einstein space,  and
examples (also non homogeneous ones) can be found in the paper 
\cite{BG}. The tangent bundle of $M^7$ splits
into the $3$-dimensional part spanned by $\eta_1, \eta_2, \eta_7$ and
its $4$-dimensional orthogonal complement. We restrict the exterior derivatives $d \eta_1, d \eta_2$ and $d \eta_7$ to this complement. In an adapted
orthonormal frame, these forms coincide with the algebraic forms $de_1, de_2$
and $de_7$. The space 
of Riemannian Killing spinors 
\bdm
\nabla^g_X \Psi \ = \ \frac{1}{2} \, X \cdot \Psi
\edm
is non trivial and has at least dimension three (see \cite{FKath}). 
Moreover, the proof of this fact shows that all the Riemannian Killing spinors
are sections in the subbundle corresponding to the $\SU(2)$-representation
$\Delta_7^0$. Now we apply Proposition \ref{3Sas} and we obtain the following 
result. 
\begin{thm}\label{thm-3-Sasa}
Let $M^7$ be $3$-Sasakian manifold in dimension seven 
and fix a Riemannian Killing spinor $\Psi$. Then there exists a 
$7$-dimensional family of torsion forms $\T$ and flux forms $\F$ 
defined by the contact structures such that
\bdm
\nabla^g_X \Psi \, + \, \frac{1}{4} \cdot (X \haken \T) \cdot
\Psi \, + \, p \cdot ( X \haken \F ) \cdot \Psi \,  + \, 
q \cdot (X \wedge \F) \cdot \Psi \ = \ 0 \ .
\edm
The condition $\F \cdot \Psi = 0$  restricts
to a subfamily of dimension six. If $4 \, p \, - \, 3 \, q \, \neq \, 0$,
the condition $\T \cdot \Psi = 0$ defines again a $6$-dimensional subfamily. 
If  $4 \, p \, - \, 3 \, q \, = \, 0$, then $\T \cdot \Psi = 
(14/3) \cdot \Psi$ for any torsion form in the family.
Both constraints together imply that the spinor field $\Psi$ is necessarily
zero.
\end{thm}
\noindent
All $3$-Sasakian manifolds are Einstein. In Section $4$, we will 
generalize this family of solutions. In particular, we construct 
homogeneous solutions on certain non-Einstein
manifolds. The qualitative behavior of these solutions
does not differ from the torsion and flux forms in the $3$-Sasakian
case, but the metric is allowed to depend on several parameters, and is
hence more flexible.
%
%
\section{Killing Spinors with $4$-Fluxes on Nearly Parallel 
$\G_2$-Manifolds}
%
\noindent
Fix a spinor $\Psi \in \Delta_7$ in the $7$-dimensional spin representation,
and consider the corresponding $3$-form $\omega^3 \in \Lambda^3(\R^7)$
defined by the formula
\bdm
\omega^3(X,Y,Z) \ := \ - \, \big(X \cdot Y \cdot Z \cdot \Psi \, , \, 
\Psi \big) .
\edm
The $3$-form acts on the spinor by 
$\omega^3 \cdot \Psi = - \, 7 \cdot \Psi$ (see \cite{FKMS}). The pair
$(\omega^3 \, , \, *\, \omega^3)$ generates a $2$-dimensional parameter
space.
\begin{prop}
 The equation
\bdm
\Big(X \, + \, \frac{r}{4} \cdot (X \haken \omega^3) \, 
+ \, p \cdot (X \haken * \, \omega^3) \, + \, q \cdot 
( X \wedge * \, \omega^3) \Big) \cdot \Psi \ = \ 0 \ 
\edm
holds for all vectors $X \in \R^7$ if and only if $\ 16 \, p = - \, 4 \, 
+ \, 12 \, q \, - \, 3 \, r$. 
\end{prop} 
\begin{proof}
A direct computation using the matrices of the spin representation
yields the result.
\end{proof}
\begin{cor} Fix parameters $(p,q) \in \R^2$. 
There a $1$-parameter family of admissible pairs, namely
\bdm
\T \ = \ \Big[\frac{12 q - 16 p}{3 f} \, - \, \frac{4}{3}\Big] \cdot 
\omega^3 \quad \mathrm{and} \quad
\F \ = \ f \cdot (* \, \omega^3).
\edm
If $4 \, p \, - \, 3 \, q = 0$, then the torsion form does not
depend on the flux form,
\bdm
\T \ = \ - \, \frac{4}{3} \cdot \omega^3 \quad \mathrm{and} \quad
\F \ = \ f \cdot (* \, \omega^3).
\edm
\end{cor}
\noindent
Consider a simply connected, nearly parallel $\G_2$-manifold $M^7$. It is an
Einstein space and we normalize the metric by the condition that the
scalar curvature equals $168$. There exists a Riemannian Killing
spinor (see \cite{FKMS})
\bdm
\nabla^g_X \Psi \ = \ X \cdot \Psi \ .
\edm
The triple $(\T, \, \F, \, \Psi)$ defined above is a solution
of the Killing equation
\bdm
\nabla^g_X \Psi \, + \, \frac{1}{4} \cdot (X \haken \T) \cdot
\Psi \, + \, p \cdot ( X \haken \F ) \cdot \Psi \,  + \, 
q \cdot (X \wedge \F) \cdot \Psi \ = \ 0 \ .
\edm
In particular, nearly parallel $\G_2$-manifolds admit a torsion
form $\T$ and a flux form $\F$ such that its Riemannian
Killing spinor is a Killing spinor with respect to the pair $(\T, \, \F)$.
The $1$-parameter family has been computed in the Corollary.
%
\section{Killing Spinors with $4$-Fluxes on a Aloff-Wallach space}
%
\noindent
The goal of this section is to construct on the Aloff-Wallach manifold 
$N(1,1) =\SU(3)/S^1$  a two-parameter family of metrics $g=g_{s,y}$ 
that admits, for every
$g_{s,y}$, a large family of torsion and flux forms making a fixed
spinor parallel. 
We use the computations available in \cite[p.109 ff]{BFGK}, 
which we hence shall not reproduce here. Consider the embedding
$S^1\ra \SU(3)$ given by $e^{i\theta}\mapsto \diag(e^{i\theta}, e^{i\theta},
e^{-2i\theta})$. The Lie algebra $\su(3)$ splits into $\su(3)=\m+\R$,
where $\R$ denotes the Lie algebra of $S^1$ deduced from the given embedding.
The space $\m$ has a preferred direction, namely the subspace
$\m_0$ generated by the matric $L:=\diag(3i, -3i,0)$. Let $E_{ij}\, (i<j)$ 
be the matrix with $1$ at the place $(i,j)$ and zero 
elsewhere, and define $A_{ij}=E_{ij}-E_{ji},\,\tilde{A}_{ij}=i(E_{ij}+E_{ji})$.
We set $\m_1:=\mathrm{Lin}\{A_{12}, \tilde{A}_{12}\}$,  
$\m_2:=\mathrm{Lin}\{A_{13}, \tilde{A}_{13}\}$ and 
$\m_3:=\mathrm{Lin}\{A_{23}, \tilde{A}_{23}\}$. The sum
$\m_1\op\m_2\op\m_3$ is an algebraic
complement of $\m_0$ inside $\m$, and in fact all spaces $\m_i$ are
pairwise perpendicular  with respect to the Killing form 
$B(X,Y):=-\mathrm{Re}(\tr XY)/2$. Hence, the following formula
\bdm
\g_{s,y} \ :=\ \frac{1}{s^2}\, B\big|_{\m_0}\, + \, B\big|_{\m_1}\,+\,
 \frac{1}{y}B\big|_{\m_2}\,+\,\frac{1}{y}B\big|_{\m_3}
\edm
defines a two-parameter family of metrics on $N(1,1):=\SU(3)/S^1$. It 
is a subfamily
of the family considered in \cite[p.109 ff]{BFGK}; in particular, 
($s=1,y=2$) corresponds to the $3$-Sasakian metric that has three
Riemannian Killing spinors with Killing number $1/2$, and ($s=1,y=2/5$) is the 
Einstein metric with one Killing spinor with Killing number $-3/10$
(see \cite[Thm 12, p.116]{BFGK}). An orthonormal basis of $\m$ is given by
\bdm
X_1\, =\,A_{12},\ X_2\, =\, \tilde{A}_{12},\
X_3\, =\,\sqrt{y} A_{13},\ X_4\, =\, \sqrt{y}\tilde{A}_{13},\
X_5\, =\,\sqrt{y}A_{23},\ X_6\, =\, \sqrt{y}\tilde{A}_{23},
\edm
and $X_7 = s\cdot L/3$.
The isotropy representation $\Ad(\theta)$ leaves the vectors
$X_1,X_2$ and $X_7$ invariant, and acts as a rotation by
$3\theta$ in the $(X_3,X_4)$-plane and in the $(X_5,X_6)$-plane. 
We use the realization of the
$8$-dimensional $\Spin(7)$-representation $\Delta_7$ as given in Section $2$. 
One then checks that $\Psi_3,\Psi_4,\Psi_5$ and $\Psi_6$
are fixed under the lift $\tilde{\Ad}(\theta)$ of the isotropy representation 
to $\Spin(7)$. Thus, they define constant sections in the spinor bundle
$S=\SU(3)\x_{\tilde{\Ad}}\Delta_7$.  The Levi-Civita connection
of $N(1,1)$ is described by a map $\Lambda:\ \m\mapsto\so(7)$, whose lift
$\tilde{\Lambda}:\ \m\mapsto\spin(7)$  can be found either in 
\cite[p.112]{BFGK} or in \cite{AgFr}.\\

\noindent
In order to define a global form on $N(1,1)$, an algebraic  form on $\m$ needs
to be invariant under the isotropy representation. It turns out that there 
are precisely $13$ isotropy invariant $3$-forms on $\m$, hence the most 
general $3$-form we can consider is 
\begin{eqnarray*}
\T&:=& \alpha_3\,(X_{135}+X_{146})+\beta_3\,(X_{235}+X_{246})+ 
\gamma_3\,(X_{357}+X_{467})+ \delta_3\,(X_{145}-X_{136})\\
&+& \eps_3\,(X_{245}-X_{236})+\xi_3\,( X_{457}-X_{367})+
\mu_3\, X_{127}+\nu_3\, X_{347}\\
&+& \lambda_3\, X_{567}+\eta_3\, X_{134} +\omega_3\, X_{234}
+\pi_3\, X_{156}+\vrho_3\, X_{256}.
\end{eqnarray*}
For notational convenience, we shall
write $X_{ijk}$ for $X_i\wedge X_j\wedge X_k$, and similarly for forms of
any degree. By Hodge duality, the Ansatz for a $4$-form is 
\begin{eqnarray*}
\F& :=& \alpha_4\, X_{1234}+\beta_4\,X_{1256}+\gamma_4\,X_{3456}+\delta_4
\,X_{1347}+ \eps_4\,X_{1567}+\xi_4\,X_{2347}+\mu_4\,X_{2567}\\
&+&\nu_4\,(X_{1235}+X_{1246})+ \lambda_4\,(X_{1357}+X_{1467})
+\eta_4\,(X_{1245}-X_{1236}) \\ 
&+&\omega_4\,(X_{1457}-X_{1367})
+ \pi_4\,(X_{2457}-X_{2367})+\vrho_4\,(X_{2357}+X_{2467}).
\end{eqnarray*}
In particular, the parameter space $\Omega$ of pairs $(\T,\F)$ of
possible $3$- and $4$-forms has now six  dimensions more
than in the $3$-Sasakian case.
Notice that for the $3$-form $\T$, $X\haken\T=-(X\cdot\T+\T\cdot X)/2$,
whereas the $4$-form $\F$ satisfies
\bdm
X\haken\F\,=\,- \frac{1}{2}(X\cdot\F-\F\cdot X),\quad
X\wedge\F\,=\, \frac{1}{2}(X\cdot\F+\F\cdot X).
\edm
\begin{thm}\label{thm-A-W}
For every metric $g_{s,y}$ on $N(1,1)$ and pair $(p,q)\in\R^2$, 
there exists a $10$-dimensional affine space $\Omega'$ of  forms 
$(\T,\F)$ such
that the spinor field $\Psi_3$ satisfies the Killing spinor equation
\bdm
\nabla_X\Psi :=\nabla^g_X\Psi+\frac{1}{4}(X\haken \T)\Psi +
p\,(X\haken \F)\Psi+q\,(X\wedge \F)\Psi = 0. 
\edm
Furthermore, the additional condition $\F\cdot\Psi_3=0$ singles out a
$9$-dimensional affine subspace of $\Omega'$. 
For $4 \, p \, - \, 3 \, q \, \neq \, 0$, 
the set of forms inside $\Omega'$ satisfying $\T\cdot\Psi_3=0$  
is again a $9$-dimensional affine subspace, but its intersection
with forms such that $\F\cdot\Psi_3=0$ is empty. For $4 \, p \, - \, 3 \, q 
\, = \, 0$, there are no
$3$-forms in $\Omega'$ such that $\T\cdot\Psi=0$.
\end{thm}
\begin{proof}
Evaluating the Killing spinor equation  in all directions 
$X=X_1,\ldots,X_7$, one observes that of the resulting seven 
$8$-dimensional spinorial equations, half is trivial, hence the linear
system in $\alpha_3,\ldots,\vrho_3,\alpha_4,\ldots,\vrho_4$ to be solved
consists only of $7\x4=28$ equations (with $4$ parameters $s,y,p,q$).
This system turns out to be highly redundant. In order to state its general 
solution, we decided to express it as functions of the parameters of $\F$. 
The $10$ coefficients $\alpha_4,\,\beta_4,\,\gamma_4,\,\eps_4,\,\xi_4,\,
\nu_4,\,\eta_4,\,\omega_4,\,\pi_4,\,\vrho_4$ can be chosen freely, the three
remaining ones are given by
\bdm\tag{$*$}
\delta_4=\eps_4-2\eta_4,\quad \mu_4=\xi_4+2\nu_4,\quad \lambda_4=\pi_4.
\edm
The coefficients of $\T$ are then expressed as functions of the coefficients 
of $\F$, hence yielding $13$ formulas. These are of two types:  the
first set is independent of the metric and relatively simple,
\bdm
\beta_3=4(p+q)\pi_4,\quad \gamma_3=-4(p+q)\nu_4,\quad \delta_3=-4(p+q)\pi_4,
\quad \xi_3=-4(p+q)\eta_4,
\edm\bdm
\eta_3=4(p+q)(\xi_4+2\nu_4), \ \omega_3=-4(p+q)\eps_4,\ 
\pi_3=4(p+q)\xi_4,\ \vrho_3=-4(p+q)(\eps_4-2\eta_4).
\edm
The second set of formulas is more complicated and, in particular,
dependent on the metric parameters $s,y$,
\begin{eqnarray*}
\alpha_3 &= & -\frac{1}{3s}\left[-2-6s+3ys+4s^2-2y-4ps(-\alpha_4+\beta_4+\gamma_4+2\omega_4-\vrho_4)+12qs\,\vrho_4\right],\\
\eps_3 & =& -\frac{1}{3s}\left[2-6s+3ys-4s^2+2y-4ps(\alpha_4-\beta_4-\gamma_4+\omega_4-2\vrho_4)-12qs\,\omega_4\right],\\
\mu_3 & =& +\frac{2}{3s}\left[-1-4s^2+2y+2ps(\alpha_4-\beta_4+2\gamma_4-2\omega_4-2\vrho_4)+6 qs\,\gamma_4\right],\\
\nu_3 & =& +\frac{1}{3s}\left[4-8s^2+y-4ps(-\alpha_4-2\beta_4+\gamma_4+2\omega_4+2\vrho_4)+12qs\,\beta_4\right],\\
\lambda_3 & =& -\frac{1}{3s}\left[4-8s^2+y-4ps(2\alpha_4+\beta_4+\gamma_4+2\omega_4+2\vrho_4)-12qs\,\alpha_4\right]\,.
\end{eqnarray*}
This shows the main part of the Theorem.
The equation $\F\cdot\Psi_3=0$ yields for the coefficients of $\F$
four conditions; three of them coincide with the equations $(*)$, whilst
the last one is the linear equation
\bdm
-\alpha_4+\beta_4+\gamma_4+2\omega_4+2\vrho_4=0\,.
\edm
Surprisingly, none of the parameters $s,\,y,\,p,\,q$ occurs. The 
constraint $\T\cdot\Psi_3=0$ gives only one condition, namely,
\bdm
s(6q-8p)(-\alpha_4+\beta_4+\gamma_4+2\omega_4+2\vrho_4)= 1+y+4s^2\,. 
\edm
Since $1+y+4s^2>0$, all remaining claims follow.
\end{proof}
%
\section{Solutions for the special $(p,q)$-coupling}\noindent
%
The coupling $4 \, p \, - 3 \, q \, = \, 0$ between the different
parts involving the flux term of the Killing equation
plays a special
role (see Theorem \ref{thm-3-Sasa} and Theorem \ref{thm-A-W}). Let us 
discuss the solutions in this case in more detail. The Killing equation
reads as ($n = 7$)
\bdm
\nabla^g_X \Psi \, + \, \frac{1}{4} \cdot (X \haken \T) \cdot
\Psi \, + \, \frac{3}{4} \cdot ( X \haken \F ) \cdot \Psi \,  + \, 
(X \wedge \F) \cdot \Psi \ = \ 0 \ .
\edm
The first series of examples are nearly parallel $\G_2$-manifolds. We 
normalize the scalar curvature by the condition $\mathrm{Scal} = 168$.
Then there exists a Riemannian Killing spinor $\Psi$ corresponding to
the $\G_2$-structure $\omega^3$,
\bdm
\nabla^g_X \Psi \ = \ X \cdot \Psi \, , \quad \omega^3 \cdot \Psi \ = \ 
- \, 7 \cdot \Psi \, .
\edm
The pair $3 \cdot \T = - \, 4 \cdot \omega^3$ and 
$\F =  f \cdot (* \omega^3)$ 
together with the spinor $\Psi$ solves the equation, where 
$f \in \R^1$ is an arbitrary real parameter. The torsion form
has a geometric meaning. It defines the unique linear, metric
connection $\nabla = \nabla^g \, + \ (1/2) \cdot \T$ preserving
the nearly parallel $\G_2$-structure 
(see \cite[Example $5.2$]{FriedrichIvanov}). Moreover, the
spinor field $\Psi$ is $\nabla$-parallel and the Killing equation decouples
into
\bdm
\nabla^g_X \Psi \, + \, \frac{1}{4} \cdot (X \haken \T) \cdot
\Psi \ = \ 0  \quad \mathrm{and} \quad  3 \cdot ( X \haken \F ) 
\cdot \Psi \,  + \, 4 \cdot (X \wedge \F) \cdot \Psi \ = \ 0 \, .
\edm 
Compact nearly parallel $\G_2$-manifolds are studied, for example, 
in \cite{FKMS}.
The $1$-parameter family of flux forms associated with a nearly
parallel $\G_2$-manifold has already been investigated in 
supergravity (see \cite{BilalDS}).\\

\noindent
A larger family of solutions arises from a $7$-dimensional $3$-Sasakian
manifold $M^7$. It is an Einstein space, and the scalar curvature is normalized
automatically to $\mathrm{Scal} = 42$. There exist four Riemannian Killing
spinors (see \cite{FKath}). Let us fix one of them. In the family
of torsion and flux forms considered in Section $3$, there exists a 
$7$-dimensional affine subspace of solutions. The torsion forms
are completely determined by the flux forms ( $p = 3/4$ and $q = 1$ in the 
notation of Section $2$),
\begin{eqnarray*}
t_{12} & = & - \, 7 \, f_{272} \ = \ - \, t_{21} , \quad 
t_{17} \ = \ - \, 7 \, f_{277} \ = \ - \, t_{71} , \\
t_{27} & = & 7 \, f_{177} \ = \ t_{72} , \\
t_{11} & = & - \, \frac{2}{3} \, + \,  f \, - \, 2 \, f_{127} \, + \, 
2 \, f_{172} \, - \, 5 \, f_{271} , \\
t_{22} & = & \frac{2}{3} \, - \, f \, + \, 
2 \, f_{127} \, - \, 2 \, f_{271} \, + \, 
5 \, f_{172} , \\
t_{77} & = & \frac{2}{3} \, - \, f \, - \, 
2 \, f_{172} \, - \, 2 \, f_{271} \, - \, 
5 \, f_{127} , \\
t & = & \frac{2}{3} \, + \, 2 \, f_{127} \, - \, 
2 \, f_{172} \, - \, 2 \, f_{271} \, + \, 6 \, f_ , \\
f_{122} & = & - \, f_{177} , \quad 
f_{121} \ = \ - \, f_{277} , \quad 
f_{171} \ = \ f_{272} \ .
\end{eqnarray*}
All torsion forms in the $7$-dimensional family of solutions act 
on the spinor by the formula $\T \cdot \Psi = (14/3) \cdot \Psi$.
The equation $\F \cdot \Psi = 0$ defines a $6$-dimensional
affine subspace,
\bdm
\gamma \ := \ f \, - \, 2 \, f_{127} \, + \, 2 \, f_{172} \, + 
\, 2 \, f_{271} \ = \ 0 \, .
\edm
We compute the action of the symmetric endomorphisms $\T$ and $\F$ on the
$8$-dimensional space of spinors explicitely. In order to formulate
the result, let us introduce the following $(3 \times 3)$-matrix $\F^*$,
\bdm
\left[\ba{ccc}f-2f_{127} -  2 f_{172} - 
2 f_{271} & 
- \, 4 \, f_{177} & - \, 4 \, f_{277} \\ 
- \, 4 \, f_{177} & f +  2 f_{127} +  2 f_{172} -  
2 f_{271} & 4 \, f_{272} \\  
- \, 4 \, f_{277} & 4 \, f_{272} & f + 2  f_{127}  - 
2 f_{172} +  2 f_{271} \ea\right] . 
\edm  
We order the basis in $\Delta_7$ in such a way that $\Psi$ is 
the last element in the basis. Then the traceless and symmetric 
endomorphisms $\T$ and $\F$ are given by the matrices
\bdm
\F \, = \, \left[\ba{ccc}-  f \cdot \mathrm{Id}_4 & 0 & 0 \\ 
0 & \F^* & 0 \\  
0 & 0 & \gamma \ea\right], \quad 
\T \, = \, \left[\ba{ccc} (\gamma - \frac{2}{3}) \cdot \mathrm{Id}_4 
& 0 & 0 \\ 
0 & (\gamma - \frac{2}{3}) \cdot \mathrm{Id}_3  & 0 \\  
0 & 0 & \frac{14}{3} -  7 \, \gamma \ea\right]\,  + \, 7 \, \F \, . 
\edm  
Remark that $\F^*$ is
an arbitrary symmetric $(3 \times 3)$-matrix. It acts in the
$3$-plane generated by the Riemannian Killing spinors orthogonal
to the fixed Riemannian Killing $\Psi$. Let us look at the family
of solutions  from the point of view of $\G_2$-structures.
The spinor $\Psi$ defines such a structure on $M^7$ (see Section 2).
Since it is a Riemannian Killing spinor, the $\G_2$-structure is
nearly parallel (see \cite{FKMS}). On the other side, a $3$-Sasakian
structure on $M^7$ is topologically a $\SU(2)$-reduction of the frame
bundle. Since $\SU(2) \subset \G_2 \subset \SO(7)$, any $3$-Sasakian
manifold induces a family of $\G_2$-structures. The spinor $\Psi$ 
singles out one of them. In any case, we have an underlying $\G_2$-structure
$\omega^3$ on $M^7$. In our parametrization of the family $(\T, \, \F)$ 
the case $\F \, = \, 0$ yields again the canonical torsion form of the
unique connection preserving the nearly parallel $\G_2$-structure
(see again \cite[Example 5.1]{FriedrichIvanov}). Moreover, the condition
\bdm
\T \ = \ \left[\ba{ccc} - \, \frac{2}{3} \cdot \mathrm{Id}_4 
& 0 & 0 \\ 
0 & - \, \frac{2}{3} \cdot \mathrm{Id}_3  & 0 \\  
0 & 0 & \frac{14}{3}  \ea\right] \
\edm
defines a $1$-parameter subfamily of flux forms. This is exactly
the above mentioned solution line of the nearly parallel $\G_2$-structure.
Consequently, if the nearly parallel $\G_2$-structure arises from
an underlying $3$-Sasakian geometry, we can embed the canonical
solution $(\T = \omega^3, \, \F = f \cdot (*\omega^3))$ into a larger
family of solutions. Then the Killing equation does not decouple anymore.
We have the same picture  for the solutions on $N(1,1)$. In this case,
the underlying $\G_2$-structure is not nearly parallel, but only cocalibrated
and additional parameters for the metric occur.\\

\noindent
A special coupling between the $(p,q)$-parameters in the Killing equation
with fluxes occurs in any dimension. We explain one way to understand this
effect. First of all, 
one easily verifies the following algebraic formulas concerning the action
of exterior forms of degree three and four on spinors:
\begin{eqnarray*}
\sum_{i=1}^n e_i \cdot \big( e_i \haken \T \big) =  3 \cdot \T  , 
\quad
\sum_{i=1}^n e_i \cdot \big( e_i \haken \F \big) =  4 \cdot \F , 
\quad \sum_{i=1}^n e_i \cdot \big( e_i \wedge \F \big) =  - \, 
(n \, - \, 4) \cdot \F \, . 
\end{eqnarray*} 
Contracting the equation
\bdm
\nabla^g_X \Psi \, + \, \frac{1}{4} \cdot (X \haken \T) \cdot
\Psi \, + \, p \cdot ( X \haken \F ) \cdot \Psi \,  + \, 
q \cdot (X \wedge \F) \cdot \Psi \ = \ 0 \, ,
\edm
we obtain
\bdm
D^g \Psi \, + \, \frac{3}{4} \, \T \cdot \Psi \, + \, 
\big( 4 \, p \, - \, (n-4) \, q \big) \cdot \F \cdot \Psi \ = \ 0 \, .
\edm
If $4 \, p \, - \, (n-4) \, q \, = \, 0$,  the
action of the Riemannian Dirac operator $D^g$ on the spinor
$\Psi$ depends only on the torsion form, but not on the flux
form. In this case we obtain a link between the spectrum of the
Riemannian Dirac operator and the admissible algebraic constraints
given by the torsion form. We can apply well-known estimates for the
Dirac spectrum of a Riemannian manifold in order to exclude some
of these solutions. For example, we obtain (see \cite{Fri1})
\begin{prop} 
Let $(M^n, \, g, \, \T, \, \F, \, \Psi)$ be a compact solution of the equation
\bdm
\nabla^g_X \Psi \, + \, \frac{1}{4} \cdot (X \haken \T) \cdot
\Psi \, + \, \frac{n-4}{4} \cdot ( X \haken \F ) \cdot \Psi \,  + \, 
(X \wedge \F) \cdot \Psi \ = \ 0 
\edm
with the constraint $\T \cdot \Psi \, = \, c \cdot \Psi$. Then 
the eigenvalue $c$ is bounded by the minimum $\mathrm{Scal}_0$
of the scalar curvature of the Riemannian manifold,
\bdm
c^2 \ \geq \ \frac{4 \, n}{9\, (n-1)} \cdot \mathrm{Scal}_0 \ . 
\edm
\end{prop}

\noindent
The solutions on $3$-Sasakian manifolds discussed before realize the 
lower bound,
since they come from a Riemannian Killing spinor. In case of the homogeneous
solutions on the Aloff-Wallach space $N(1,1)$ the eigenvalue $c^2$ is
strictly greater then the lower bound. 
\begin{NB} On an $8$-dimensional manifold the equation corresponding
to the special $(p,q)$-parameters simplifies, 
\bdm
\nabla^g_X \Psi \, + \, \frac{1}{4} \cdot (X \haken \T) \cdot \Psi \, + \,
\F \cdot X \cdot \Psi \ = \ 0 \, .
\edm
\end{NB}
\noindent
Our method for the construction of torsion and flux forms solving
the equation at hand applies in dimension eight, too. The key point is the
following Proposition. Its proof relies on a computer computation.
We have to solve a system of $128$ linear equations
in $126$ variable, and turns out to have sufficiently many solutions.
\begin{prop} 
Let $\Psi = \Psi^+ \, + \, \Psi^-$ be an $8$-dimensional
spinor with non trivial positive and negative part, $\Psi^{\pm} \neq 0$.
Then there exists a family depending on $25$ parameters of $3$-forms
$\T \in \Lambda^3(\R^8)$ and $4$-forms $\F \in \Lambda^4(\R^8)$ such that,
for any vector $X \in \R^8$, the following equation holds:
\bdm
\frac{1}{4} \cdot (X \haken \T) \cdot \Psi \, + \,
\F \cdot X \cdot \Psi \ = \ 0 \ .
\edm
\end{prop}
\noindent
Consider an $8$-dimensional Lie group $(\G^8, \, g)$ equipped
with a biinvariant Riemannian metric. The formula
\bdm
\T_0(X,Y,Z) \ := \ -\, g( [X\, , \, Y] \, , \, Z)
\edm
defines the canonical torsion form of the Lie group. The action
of the Levi-Civita connection on a spinor field $\Psi : \G^8 \rightarrow 
\Delta_8$ is given by the formula
\bdm
\nabla^g_X \Psi \ = \ d \Psi(X) \, + \, \frac{1}{4} \, 
(X \haken \T_0) \cdot \Psi \ .
\edm
Consequently, any constant spinor field on the Lie group admits a
family depending on $25$ parameters of torsion and flux forms solving
the equation.
    
\end{document}